\documentstyle[epsf]{article}
\input{prepictex}
\input{pictex}
\input{postpictex}
\begin{document}
\def\thefootnote{}
\def\lhead{A. Tyszka}
\def\rhead{Beckman-Quarles theorem for rational eight-space}
\newcommand{\F}{\bf F}
\newcommand{\R}{\bf R}
\newcommand{\Q}{\bf Q}
\newcommand{\wt}{\widetilde}
\hsize=126mm
\vsize=180mm
\parindent=5mm
\vskip 0.1truecm
\centerline{\bf A discrete form of the Beckman-Quarles
theorem for rational eight-space}
\vskip 0.5truecm
\centerline{\large Apoloniusz Tyszka}
\vskip 0.5truecm
\vskip 0.5truecm
\footnotetext{
{\bf AMS (1991) Subject Classification:} Primary 51M05, Secondary 05C12}
\hfill
\par
{\bf Abstract}. Let $\Q$ be the field of rationals numbers.
We prove that:
{\large (1)}
if $x,y\in {\R}^{n}$ $(n>1)$ and $|x-y|$ is constructible by means of
ruler and compass then there exists a finite set
$S_{xy}\subseteq {\R}^{n}$ containing $x$ and $y$ such that each map
from $S_{xy}$ to ${\R}^{n}$ preserving unit distance preserves
the distance between $x$ and $y$,
{\large (2)}
if $x,y \in {\Q}^{8}$ then there exists a finite set
$S_{xy} \subseteq {\Q}^{8}$ containing $x$ and $y$ such that each map
from $S_{xy}$ to ${\R}^{8}$ preserving unit distance preserves
the distance between $x$ and $y$.
\vskip 0.5truecm
Theorem 1 may be viewed as a discrete form of the classical
Beckman-Quarles theorem, which states that any map from
${\R}^{n}$ to ${\R}^{n}$ ($2\leq n<\infty $) preserving unit
distance is an isometry, see [1]-[3]. Theorem 1 was announced
in [9] and prove there in the case where $n=2$. A stronger version of
Theorem 1 can be found in [10], but we need the elementary proof
of Theorem 1 as an introduction to Theorem 2.
\vskip 0.1truecm
{\bf Theorem 1}. If $x,y\in {\R}^{n}$ ($n>1$) and $|x-y|$ is
constructible by means of ruler and compass then there exists a
finite set $S_{xy}\subseteq {\R}^{n}$ containing $x$ and $y$ such
that each map from $S_{xy}$ to ${\R}^{n}$ preserving unit
distance preserves the distance between $x$ and $y$.
\vskip 0.1truecm
{\it Proof.}
Let us denote by $D_{n}$ the set of all non-negative numbers $d$
with the following property:
\vskip 0.1truecm
\noindent
if $x,y\in {\R}^{n}$ and $|x-y|=d$ then there exists a finite set
$S_{xy}\subseteq {\R}^{n}$ such that $x,y\in S_{xy}$ and any map
$f:S_{xy}\rightarrow {\R}^{n}$ that preserves unit distance
preserves also the distance between $x$ and $y$.
\vskip 0.1truecm
Obviously $0,1\in D_{n}$. We first prove that if $d\in D_{n}$ then
$\sqrt{2+2/n}\cdot d\in D_{n}$. Assume that $d>0$,
$x,y\in {\R}^{n}$ and $|x-y|=\sqrt{2+2/n}\cdot d$.
Using the notation of Figure 1 we show that
\vskip 0.2truecm
\centerline{$S_{xy}:=
\bigcup \{S_{ab}:a,b\in \{x,y,\wt{y},p_{1},p_{2},...,p_{n},
\wt{p}_{1},\wt{p}_{2},...,\wt{p}_{n}\},|a-b|=d\}$}
\vskip 0.2truecm
\noindent
satisfies the condition of Theorem 1.
Figure 1 shows the case $n=2$, but equations below Figure 1
describe the general case $n\geq 2$; $z$ denotes
the centre of the $(n-1)$-dimensional regular simplex
$p_{1}p_{2}...p_{n}$.
\\
\\
\\
\centerline{\epsfbox{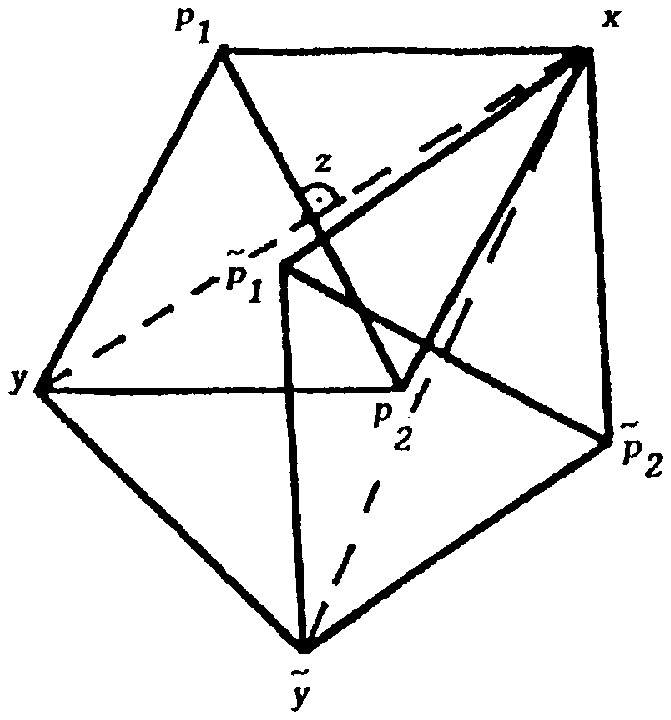}}
\centerline{Figure 1}
\centerline{$1\leq i<j\leq n$}
\centerline{$|y-\wt{y}|=d$,
$|x-p_{i}|=|y-p_{i}|=|p_{i}-p_{j}|=d=|x-\wt{p}_{i}|=
|\wt{y}-\wt{p}_{i}|=|\wt{p}_{i}-\wt{p}_{j}|$}
\centerline{
$|x-\wt{y}|=|x-y|=2\cdot |x-z|=2\cdot \sqrt{\frac{n+1}{2n}}
\cdot d=\sqrt{2+2/n}\cdot d$}
\vskip 0.5truecm
Assume that $f:S_{xy}$ $\rightarrow$ ${\R}^{n}$
preserves distance $1$. Since
\vskip 0.3truecm
\centerline{$S_{xy}\supseteq S_{y\wt{y}}
\cup
\bigcup_{i=1}^{n}S_{xp_{i}}
\cup
\bigcup_ {i=1}^{n}S_{yp_{i}}
\cup
\bigcup_{1 \leq i<j \leq n} S_{p_{i}p_{j}}$}
\vskip 0.3 truecm
\noindent
we conclude that $f$ preserves the distances between $y$ and $\wt{y}$,
$x$ and $p_{i}$ ($1\leq i \leq n$), $y$ and $p_{i}$ ($1\leq i\leq n$),
and all distances between $p_{i}$ and $p_{j}$ ($1\leq i<j\leq n$).
Hence $|f(y)-f(\wt{y})|=d$ and $|f(x)-f(y)|$ equals either $0$ or
$\sqrt{2+2/n}\cdot d$. Analogously we have that $|f(x)-f(\wt{y})|$
equals either $0$ or $\sqrt{2+2/n}\cdot d$. Thus $f(x)\neq f(y)$, so
$|f(x)-f(y)|=\sqrt{2+2/n}\cdot d$ which completes the proof that
$\sqrt{2+2/n}\cdot d\in D_{n}$.
\vskip 0.3truecm
\noindent
Therefore, if $d\in D_{n}$ then
$(2+2/n)\cdot d=\sqrt{2+2/n}\cdot (\sqrt{2+2/n}\cdot d)\in D_{n}$.
\vskip 0.3truecm
We next prove that if $x,y\in {\R}^{n}$, $d\in D_{n}$
and $|x-y|=(2/n)\cdot d$ then there exists a finite set
$Z_{xy}\subseteq {\R}^{n}$ containing $x$ and $y$
such that any map $f:Z_{xy} \rightarrow {\R}^{n}$ that preserves
unit distance satisfies $|f(x)-f(y)|\leq |x-y|$; this result is adapted
from [3]. It is obvious in the case where $n=2$,
therefore we assume that $n>2$ and $d>0$.
In Figure 2, $z$ denotes the centre of the
$(n-1)$-dimensional regular simplex $p_{1}p_{2}...p_{n}$.
Figure 2 shows the case $n=3$, but equations below Figure 2 describe
the general case where $n\geq 3$.
\\
\centerline{\epsfbox{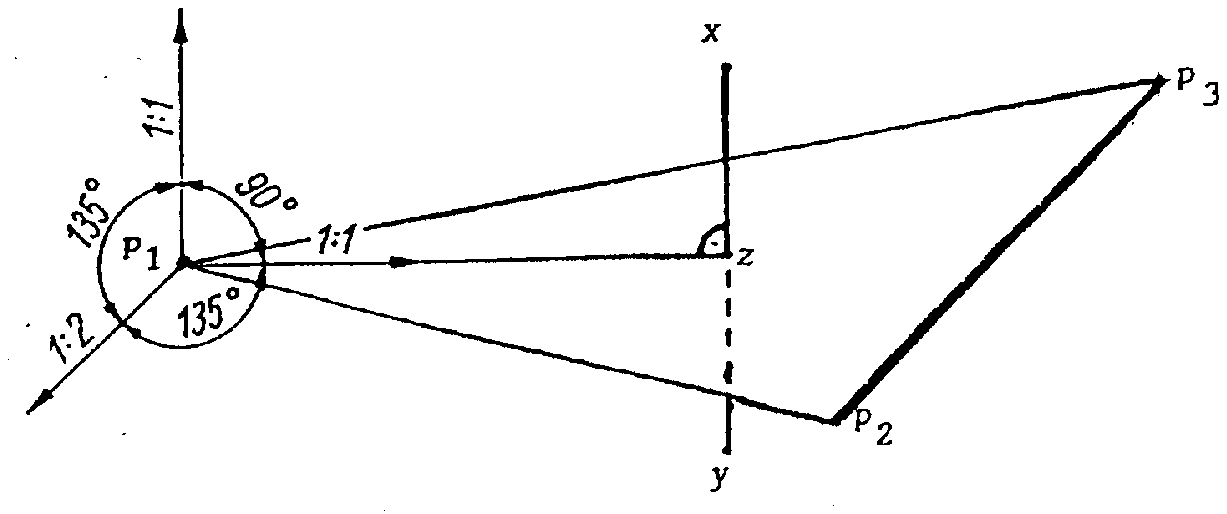}}
\centerline{Figure 2}
\vskip 0.1truecm
\centerline{$1 \leq i<j \leq n$}
\vskip 0.1truecm
\centerline{$|x-p_{i}|=|y-p_{i}|=d, \hspace{0.20cm}
|p_{i}-p_{j}|=\sqrt{2+2/n} \cdot d, \hspace{0.20cm}
|z-p_{i}|=\sqrt{1-1/n^{2}} \cdot d$}
\vskip 0.1truecm
\centerline{$|x-y|=2 \cdot |x-z|=2 \cdot
\sqrt{|x-p_{i}|^{2}-|z-p_{i}|^2}=
2 \cdot \sqrt{d^{2}-(1-1/n^{2}) \cdot d^{2}}=(2/n) \cdot d$}
\vskip 0.4truecm
\noindent
Define:
\vskip 0.2truecm
\centerline{$Z_{xy}:=\bigcup_{1\le i<j \leq n} S_{p_{i}p_{j}}
\cup
\bigcup_{i=1}^{n} S_{xp_{i}}
\cup
\bigcup_{i=1}^{n} S_{yp_{i}}$}
\vskip 0.2truecm
\noindent
If $f:Z_{xy}\rightarrow {\R}^{n}$ preserves distance $1$ then
$|f(x)-f(y)|=|x-y|=(2/n)\cdot d$ or $|f(x)-f(y)|=0$, hence
$|f(x)-f(y)|\leq|x-y|$.
\vskip 0.2truecm
\noindent
If $d\in D_{n}$, then $2\cdot d\in D_{n}$ (see Figure 3).
\\
\centerline{\epsfbox{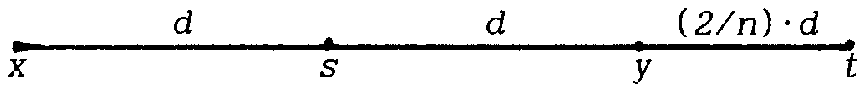}}
\centerline{Figure 3}
\centerline{$|x-y|=2 \cdot d$}
\centerline{$S_{xy}=S_{xs}\cup S_{sy}\cup Z_{yt}\cup S_{xt}$}
\vspace{0.1cm}
\\
From Figure 4 it is clear that if $d\in D_{n}$ then all distances $k\cdot d$
(where $k$ is a positive integer) belong to $D_{n}$.
\\
\centerline{\epsfbox{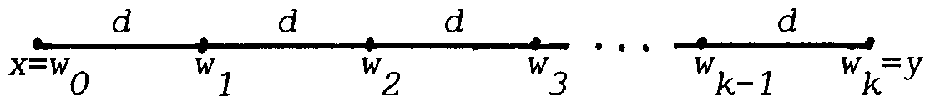}}
\centerline{Figure 4}
\centerline{$|x-y|=k\cdot d$}
\centerline{$S_{xy}=\bigcup \{S_{ab}:a,b\in \{w_{0},w_{1},...,w_{k}\},
|a-b|=d \vee |a-b|=2\cdot d\}$}
\vskip 0.2truecm
\noindent
From Figure 5 it is clear that if $d\in D_{n}$ then all distances
$d/k$ (where $k$ is a positive integer) belong to $D_{n}$.
Hence $\Q$ $\cap$ $(0,\infty)$ $\subseteq$ $D_{n}$.
\\
\centerline{\epsfbox{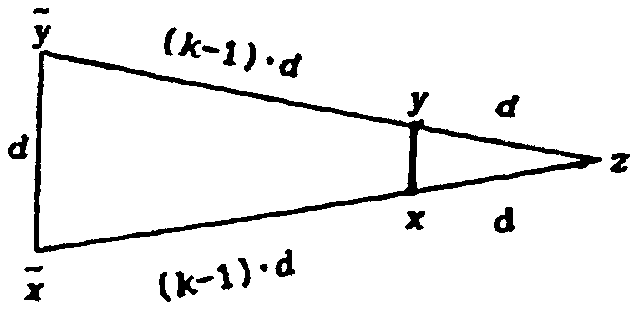}}
\centerline{Figure 5}
\centerline{$|x-y|=d/k$}
\centerline{$S_{xy}=
S_{\wt{x}\wt{y}}
\cup
S_{\wt{x}x}\cup
S_{xz}\cup
S_{\wt{x}z}
\cup
S_{\wt{y}y}
\cup S_{yz}
\cup
S_{\wt{y}z}$}
\vskip 0.5truecm
{\bf Observation}. If $x,y\in {\R}^{n}$ ($n>1$) and $\varepsilon >0$
then there exists a finite set $T_{xy}(\varepsilon) \subseteq
{\R}^{n}$ containing $x$ and $y$ such that for each map
$f:T_{xy}(\varepsilon )\rightarrow {\R}^{n}$
preserving unit distance we have
$||f(x)-f(y)|-|x-y||\leq \varepsilon $.
\\
{\it Proof}. It follows from Figure 6.
\\
\centerline{\epsfbox{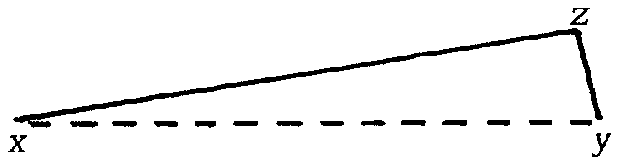}}
\centerline{Figure 6}
\centerline{$|x-z|,|z-y|\in{\Q} \cap (0,\infty)$,
\hspace{0.2cm}
$|z-y|\leq \varepsilon/2$}
\centerline{$T_{xy}(\varepsilon)=S_{xz}\cup S_{zy}$}
\vskip 0.3truecm
{\bf Note.} The above part of the proof can be found in [10].
\vskip 0.4 truecm
\noindent
If $a,b\in D_{n}$, $a>b>0$ then
$\ \sqrt{a^{2}-b^{2}} \in D_{n}$ (see Figure 7, cf.[9]).
\\
\centerline{
\beginpicture
\Large
\setcoordinatesystem units <1mm, 1mm>
\setplotarea x from 0 to 28, y from -2.5 to 19.5
\setplotsymbol({\xxpt\rm.})
\plot 14 0 14 17 0 0 28 0 14 17 /
\put {$s$} at 0 -2.5
\put {$x$} at 14 -2.5
\put {$t$} at 28 -2.5
\put {$y$} at 14 19.5
\put {$a$} at 5.5 10
\put {$a$} at 22.5 10
\put {$b$} at 7 2.5
\put {$b$} at 21 2.5
\circulararc 90 degrees from 14 2.5 center at 14 0
\put {.} at 13 1
\endpicture}
\centerline{Figure 7}
\centerline{$|x-y|$=$\sqrt{a^{2}-b^{2}}$}
\centerline{
$S_{xy}=S_{sx}\cup S_{xt}\cup S_{st}\cup S_{sy}\cup S_{ty}$}
\vskip 0.5truecm
\noindent
Hence
$\sqrt{3}\cdot a=\sqrt{(2\cdot a)^{2}-a^{2}}\in D_{n}$ and $\sqrt{2}
\cdot a=\sqrt{(\sqrt{3}\cdot a)^{2}-a^{2}}\in D_{n}$.
Therefore $\sqrt{a^{2}+b^{2}}=
\sqrt{(\sqrt{2}\cdot a)^{2}-(\sqrt{a^{2}-b^{2}})^{2}}\in D_{n}$.
\vskip 0.7truecm
In Figure 8, z denotes the centre of the
$(n-1)$-dimensional regular simplex $p_{1}p_{2}...p_{n}$, $n=2$,
but equations below Figure 8 describe the general case where $n\geq 2$.
This construction shows that if $a,b\in D_{n}$, $a>b>0$, $n\geq 2$
then $a-b\in D_{n}$, hence $a+b=2\cdot a-(a-b)\in D_{n}$.
\\
\centerline{
\beginpicture
\LARGE
\setcoordinatesystem units <1mm, 1mm>
\setplotarea x from 0 to 40, y from -3 to 40
\setplotsymbol({\xxpt\rm.})
\plot 0 0 20 36 40 0 0 0 20 15 40 0 /
\setdashes <2mm>
\plot 20 0 20 36 /
\put {$z$} at 20 -3
\put {$y$} at 22 17
\put {$x$} at 20 38
\put {$p_1$} at 0 -3
\put {$p_2$} at 40 -3
\setsolid
\circulararc 90 degrees from 20 2.5 center at 20 0
\put {.} at 19 1
\endpicture}
\centerline{Figure 8}
\centerline{$|x-y|=a-b,\hspace{0.2cm} |x-z|=a \in D_{n},\hspace{0.2cm}
|y-z|=b \in D_{n}$}
\centerline{$|p_{i}-p_{j}|=\sqrt{2+2/n}\in D_{n},\hspace{0.2cm}
|z-p_{i}|=\sqrt{1^{2}-(1/n)^{2}}\in D_{n},\hspace{0.2cm}
1 \leq i<j \leq n$}
\centerline{$|x-p_{1}|=\sqrt{|x-z|^{2}+|z-p_{1}|^{2}}=...
=|x-p_{n}|=\sqrt{|x-z|^{2}+|z-p_{n}|^{2}}\in D_{n}$}
\centerline{$|y-p_{1}|=\sqrt{|y-z|^{2}+|z-p_{1}|^{2}}=...
=|y-p_{n}|=\sqrt{|y-z|^{2}+|z-p_{n}|^{2}}\in D_{n}$}
\centerline{
$S_{xy}=\bigcup_{1\le i<j \le n}S_{p_{i}p_{j}}
\cup
\bigcup_{i=1}^{n}S_{xp_{i}}
\cup
\bigcup_{i=1}^{n}S_{yp_{i}}
\cup
T_{xy}(b)$}
\vskip 0.1truecm
In order to prove that $D_{n}\backslash \{0\}$ is a multiplicative
group it remains to observe that if positive
$a,b,c\in D_{n}$, then $\frac{a\cdot b}{c} \in D_{n}$
(see Figure 9, cf.[9]).
\\
\centerline{\epsfbox{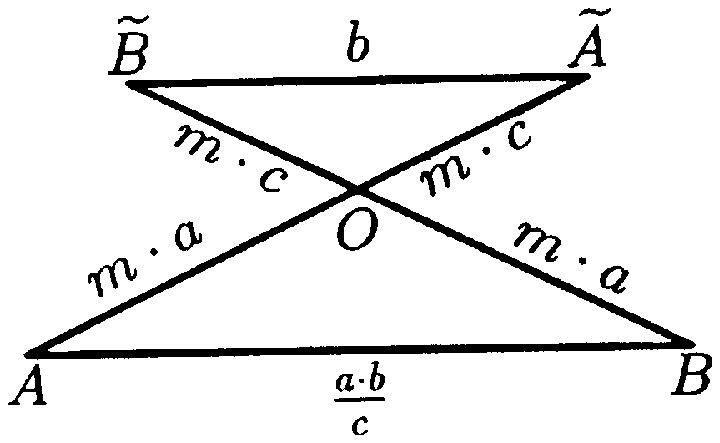}}
\centerline{Figure 9}
\centerline{$m$ is a positive integer}
\centerline{$b<2\cdot m\cdot c$}
\centerline{$S_{AB}=S_{OA}
\cup
S_{OB}\cup
S_{O\wt{A}}
\cup
S_{O\wt{B}}
\cup
S_{A\wt{A}}
\cup
S_{B\wt{B}}
\cup
S_{\wt{A}\wt{B}}$}
\vskip 0.6 truecm
If $a\in D_{n}$, $a>1$, then $\sqrt{a}=\frac{1}{2}\cdot
\sqrt{(a+1)^{2}-(a-1)^{2}} \in D_{n}$;
if $a\in D_{n}$, $0<a<1$,
then $\sqrt{a}=1/\sqrt{\frac{1}{a}}\in D_{n}$.
Thus $D_{n}$ contains all non-negative real numbers contained
in the real quadratic closure of $\Q$. This completes the proof.
\vskip 0.1truecm
\noindent
{\bf Remark 1.} Let $\F \subseteq \R$ is a euclidean field, i.e.
$ \forall x \in \F \exists {\rm y} \in \F$ ($x=y^{2} \vee x=-y^{2}$)
(cf. [6]). Our proof of Theorem 1 gives that if
$x,y \in {\F}^{n}$ ($n>1$) and $|x-y|$ is constructible by
means of ruler and compass then there exists a finite set
$S_{xy} \subseteq {\F}^{n}$ containing $x$ and $y$ such that each
map from $S_{xy}$ to ${\R}^{n}$ preserving unit distance
preserves the distance between $x$ and $y$.
\\
\\
{\bf Theorem 2.} If $x,y \in {\Q}^{8}$ then there exists a finite set
$S_{xy} \subseteq {\Q}^{8}$
containing $x$ and $y$ such that each map from $S_{xy}$ to ${\R}^{8}$
preserving unit distance preserves the distance between $x$ and $y$.
\\
\\
{\it Proof.} Denote by $R_{8}$ the set of all $d \geq 0$ with
the following property:
\\
\\
if $x,y \in {\Q}^{8}$ and $|x-y|=d$ then there exists a finite set
$S_{xy} \subseteq {\Q}^{8}$ such
that $x,y \in S_{xy}$ and any map $f:S_{xy} \rightarrow {\R}^{8}$
that preserves unit distance preserves also the distance between $x$
and $y$.
\\
\\
Obviously $0,1 \in R_{8}$. We need to prove that if
$x \neq y \in {\Q}^{8}$ then $|x-y| \in R_{8}$.
We show that configurations from
Figures 1-5 and 7 (see the proof of Theorem 1) exist in
${\Q}^{8}$. We start from simple lemmas.
\\
\\
{\bf Lemma 1} (see [11]). If $A$ and $B$ are two different points of
${\Q}^{n}$ then the reflection of ${\Q}^{n}$ with respect to the
hyperplane which is the perpendicular bisector of the segment $AB$,
is a rational transformation (that is, takes rational points to
rational points).
\\
\vskip 0.05truecm
\noindent
{\bf Lemma 2} (in the real case cf.[2] p.173 and [10]).
If $A,B,\wt{A},\wt{B} \in {\Q}^{8}$ and $|AB|=|\wt{A}\wt{B}|$
then there exists an isometry
$I:{\Q}^{8} \rightarrow {\Q}^{8}$ satisfying $I(A)=\wt{A}$
and $I(B)=\wt{B}$.
\\
{\it Proof.} If $A=\wt{A}$ and $B=\wt{B}$ then $I=id({\Q}^{8})$.
If $A=\wt{A}$ and $B \neq \wt{B}$ then the reflection of
${\Q}^{8}$ with respect to the hyperplane which is the perpendicular
bisector of the segment $B\wt{B}$, satisfies the condition of Lemma 2
in virtue of Lemma 1. Assume that $A \neq \wt{A}$.
Let $I_{1}:{\Q}^{8} \rightarrow {\Q}^{8}$ denote the reflection
of ${\Q}^{8}$ with respect to the hyperplane which is the perpendicular
bisector of the segment $A\wt{A}$. If
$I_{1}(B)=\wt{B}$ then the proof is complete. In the opposite case let
$B_{1}=I_{1}(B)$,
$B_{1} \in {\Q}^{8}$ according to Lemma 1. Let
$I_{2}:{\Q}^{8} \rightarrow {\Q}^{8}$
denote the
reflection of ${\Q}^{8}$ with respect to the hyperplane which is the
perpendicular bisector of the segment $B_{1}\wt{B}$.
Since $|\wt{A}B_{1}|=|I_{1}(A)I_{1}(B)|=|AB|=|\wt{A}\wt{B}|$ we conclude
that $I_{2}(\wt{A})=\wt{A}$. Therefore $I=I_{2} \circ I_{1}$ satisfies
the condition of Lemma 2.
\vskip 0.05truecm
{\bf Corollary.} Lemma 2 ensures that if some configuration from
Figures 1-5 and 7 exists in ${\Q}^{8}$ for a fixed $x,y \in {\Q}^{8}$,
then this configuration exists for any $x,y \in {\Q}^{8}$ with the
same $|x-y|$.
\\
\\
{\bf Lemma 3.} LAGRANGE'S FOUR SQUARE THEOREM. Every non-negative
integer is the sum of four squares of integers, and therefore
every non-negative rational is the sum of four squares of
rationals, see [8].
\\
\\
{\bf Lemma 4.} If $a,b$ are positive rationals and $b<2a$ then
there exists a triangle in ${\Q}^{8}$ with sides $b, a, a$.
\\
\\
{\it Proof.} Let $a^{2}-{(b/2)}^{2}=k^{2}+l^{2}+m^{2}+n^{2}$
where $k,l,m,n$ are rational according to Lemma 3.
Then the triangle
\\
\centerline
{$[-b/2,0,0,0,0,0,0,0]$ $[b/2,0,0,0,0,0,0,0]$ $[0,k,l,m,n,0,0,0]$}
\\
has sides $b,a,a$.
\\
\\
Now we turn to the main part of the proof. Rational coordinates
of the following configuration are taken from [11].
\\
\\
\noindent
$x=[0,0,0,0,0,0,0,0]$
\\
\\
$y=(3/8) \cdot[-3,0,0,0,1,1,1,-2]$
\hspace{0.80cm}
$\wt{y}=(1/6) \cdot [-8,1,1,3,1,0,-1,-2]$
\\
\\
$p_{1}=[-1,0,0,0,0,0,0,0]
\hspace{2.10cm}
p_{2}=(1/2) \cdot [-1,1,0,0,0,0,1,-1]$
\\
\\
$p_{3}=(1/2) \cdot [-1,-1,0,0,0,0,1,-1]$
\hspace{0.50cm}
$p_{4}=(1/2) \cdot [-1,0,1,0,0,1,0,-1]$
\\
\\
$p_{5}=(1/2) \cdot [-1,0,-1,0,0,1,0,-1]$
\hspace{0.50cm}
$p_{6}=(1/2) \cdot [-1,0,0,1,1,0, 0,-1]$
\\
\\
\hspace{0.50cm}
$p_{7}=(1/2) \cdot [-1,0,0,-1,1,0,0,-1]$
\hspace{0.50cm}
$p_{8}=(1/2) \cdot [-1,0,0,0,1,1,1,0]$
\\
\\
\\Let $I:{\Q}^{8} \rightarrow {\Q}^{8}$ denote the reflection with respect
to the hyperplane which is the perpendicular bisector of the segment
$y\wt{y}$.
By Lemma 1 we have $\wt{p}_{i}=I(p_{i}) \in {\Q}^{8}$ ($ 1 \leq i\leq 8 $).
It is easy to check that points $x,y,\wt{y},p_{i},\wt{p}_{i}$
($1\leq i \leq 8 $) form the configuration
from Figure 1 for $d=1$. The Corollary ensures that
$3/2=\sqrt{2+2/8} \cdot d=|x-y| \in R_{8}$.
\\
\vskip 0.1truecm
Points $(3/2)x, (3/2)y, (3/2)\wt{y}, (3/2)p_{i}, (3/2)\wt{p_{i}}$
($1 \leq i \leq 8$) form the configuration from Figure 1
for $d=\sqrt{2+2/8}=3/2$. The Corollary ensures that
$2+1/4=\sqrt{2+2/8} \cdot d=|(3/2) \cdot x - (3/2) \cdot y| \in R_{8}$. 
\\
\\
The following points:
\vskip 0.1truecm
$p_1=[-3/2,0,0,0,0,0,0,0]$
\vskip 0.1truecm
$p_{2}=[-3/4,3/4,0,0,0,0,3/4,-3/4]$
\vskip 0.1truecm
$p_{3}=[-3/4,-3/4,0,0,0,0,3/4,-3/4]$
\vskip 0.1truecm
$p_{4}=[-3/4,0,3/4,0,0,3/4,0,-3/4]$
\vskip 0.1truecm
$p_{5}=[-3/4,0,-3/4,0,0,3/4,0,-3/4]$
\vskip 0.1truecm
$p_{6}=[-3/4,0,0,3/4,3/4,0,0,-3/4]$
\vskip 0.1truecm
$p_{7}=[-3/4,0,0,-3/4,3/4,0,0,-3/4]$
\vskip 0.1truecm
$p_{8}=[-3/4,0,0,0,3/4,3/4,3/4,0]$
\vskip 0.1truecm
$x=[-3/4,0,0,0,1/4,1/4,1/4,-1/2]$
\vskip 0.1truecm
$y=[-15/16,0,0,0,5/16,5/16,5/16,-5/8]$
\vskip 0.1truecm
\noindent
form the configuration from Figure 2 for $d=1$. Therefore, in virtue of
Corollary if $x,y \in {\Q}^{8}$ and $|x-y|=(2/8)\cdot d = 1/4$,
then there exists a finite set $Z_{xy} \subseteq {\Q}^{8}$
containing $x$ and $y$ such that any map
$f:Z_{xy} \rightarrow {\R}^{8}$ that preserves unit distance
satisfies $|f(x)-f(y)| \leq |x-y|$.
\\
\\
As in the proof of Theorem 1 we can
prove that $2 \in R_{8}$ and all integer distances belong to
$R_{8}$. In the same way using the Corollary we can prove that
all rational distances belong to $R_{8}$, because by Lemma 4
there exists a triangle in ${\Q}^{8}$ with sides
$d,k \cdot d, k \cdot d$ ($d,k$ are positive integers, see Figure 5).
\\
\\
Finally, we prove that $|x-y| \in R_{8}$ for arbitrary
$x \neq y \in {\Q}^{8}$. It is
obvious if $|x-y|=1/2$ because $1/2$ is rational.
Let us assume that
$|x-y| \neq 1/2$. We have:
${|x-y|}^{2}$ $\in \Q \cap$ ($0,\infty$). Let
${|x-y|}^{2}=k^{2}+l^{2}+m^{2}+n^{2}$
where $k,l,m,n$ are rationals according to Lemma 3.
\\
\\
The following points:
\\
\\
$s=[-||x-y|^2-1/4|,0,0,0,0,0,0,0]$
\hspace{0.1cm},\hspace{0.1cm}
$x=[0,0,0,0,0,0,0,0]$
\\
\\
$t=[||x-y|^2-1/4|,0,0,0,0,0,0,0]$
\hspace{0.1cm},\hspace{0.1cm}
$y=[0,k,l,m,n,0,0,0]$
\\
\\
\noindent
form the configuration from Figure 7 for
$a=|{|x-y|}^{2}+1/4|$ $\in \Q \cap$ ($0,\infty$) $\subseteq R_{8}$
and
$b=|{|x-y|}^{2}-1/4|$ $\in \Q \cap$ ($0,\infty$) $\subseteq R_{8}$.
The Corollary ensures that
$|x-y|=\sqrt{a^{2}-b^{2}}\in R_{8}$.
This completes the proof of Theorem 2.
\\
\\
\noindent
{\bf Remark 2.} Theorem 2 implies that any map
$f:{\Q}^{8}\rightarrow {\Q}^{8}$
which preserves unit distance is an isometry.
\\
{\bf Remark 3.}
It is known that the injection of
${\Q}^{n}$ ($n \geq 5$)
which preserves the distances $d$ and $d/2$ ($d$ is positive and rational)
is an isometry, see [5]. The general result from [7] implies that any map
$f:{\Q}^{n} \rightarrow {\Q}^{n}$ ($n \geq 5$)
which preserves the distances $1$ and $4$ is an isometry.
On the other hand, from [4] (for $n=1,2$) and [5] (for $n=3,4$)
it may be concluded that there exist bijections of
${\Q}^{n}$ ($n=1,2,3,4$) which preserve all distances belonging to
$\{k/2: k=1,2,3,...\}$ and which are not isometries.
\\
{\bf Remark 4}. J. Zaks informed (private communication, May 2000)
the author that he proved the following:
\\
1. (cf. Remark 3): Let $k$ be any integer, $k \ge 2$;
every mapping from ${\Q}^n$ to ${\Q}^n$, $n \ge 5$,
which preserves the distances $1$ and $k$ - is an isometry.
\\
2. Theorem 2 holds for all even $n$ of the form $n=4t(t+1)$,
$t \ge 2$, as well as for all odd values of $n$ which are a
perfect square greater than $1$, $n=x^2$, and which, in addition,
are of the form $n=2y^2-1$.
The construction is a modified version of the proof of Theorem 2.

\vskip 1 truecm
\centerline{{\bf References}}
\par
\begin{enumerate}
\item F. S. Beckman and D. A. Quarles Jr.,
On isometries of euclidean spaces,
{\it Proc. Amer. Math. Soc.}, 4 (1953), 810-815.
\item W. Benz,
{\it Real geometries}, BI Wissenschaftsverlag, Mannheim, 1994.
\item U. Everling, Solution of the isometry problem stated
by K. Ciesielski, {\it Math. Intelligencer} 10 (1988), No.4, p.47.
\item B. Farrahi, On distance preserning transformations of Euclidean-like
planes over the rational field, {\it Aequationes Math.} 14 (1976), 473-483.
\item B. Farrahi, A characterization of isometries of rational
euclidean spaces, {\it J. Geom.} 12 (1979), 65-68.
\item M. Hazewinkel, {\it Encyclopaedia of mathematics}, Kluwer
Academic Publishers, Dordrecht 1995.
\item H. Lenz, Der Satz von Beckman-Quarles im rationalen Raum,
{\it Arch. Math.} (Basel) 49 (1987), 106-113.
\item L. J. Mordell, {\it Diophantine equations}, Academic Press,
London - New York, 1969.
\item A. Tyszka, A discrete form of the Beckman-Quarles theorem,
{\it Amer. Math. Monthly} 104 (1997), 757-761.
\item A. Tyszka, Discrete versions of the Beckman-Quarles theorem,
{\it Aequationes Math.} 59 (2000), 124-133.
\item J. Zaks, On the chromatic number of some rational spaces,
{\it Ars. Combin.} 33 (1992), 253-256.
\end{enumerate}
\vskip 0.1truecm
{\it Technical Faculty}
\\
{\it Hugo Ko{\l}{\l}\c{a}taj University}
\\
{\it Balicka 104, PL-30-149 Krak\'ow, Poland}
\\
{\it rttyszka@cyf-kr.edu.pl}
\\
{\it http://www.cyf-kr.edu.pl/\symbol{126}rttyszka}
\end{document}